\documentclass[12pt]{amsart}

\usepackage{graphics}
\usepackage{amsmath}
\usepackage{amsfonts}
\usepackage{amssymb}
\usepackage{amscd}

\font\Bbb=msbm10

\def\C{\hbox{\Bbb C}}

\def\Z{\hbox{\Bbb Z}}
\def\N{\hbox{\Bbb N}}
\def\e{{\varepsilon}}
\def\pd#1{ \partial_{#1} }

\newtheorem{theorem}{\bf Theorem}
\newtheorem{lemma}[theorem]{\bf Lemma}
\newtheorem{remark}[theorem]{\bf Remark}
\newtheorem{example}[theorem]{\bf Example}
\newtheorem{proposition}[theorem]{\bf Proposition}
\newtheorem{corollary}[theorem]{\bf Corollary}
\newtheorem{definition}[theorem]{\bf Definition}

\numberwithin{equation}{section} \numberwithin{theorem}{section}

\title{Bases in the solution space of the Mellin system}

\author{Alicia Dickenstein}
\address{Dto. de Matem\'atica, FCEN,
Universidad de Buenos Aires, \newline \indent (1428) Buenos Aires,
Argentina.} \email{alidick@dm.uba.ar}
\thanks{Alicia Dickenstein was partially supported  by
UBACYT X042, CONICET PIP 5617 and ANPCYT PICT 20569, Argentina.}

\author{Timur Sadykov}
\address{Department of Mathematics and Computer Science, \newline
\indent Krasnoyarsk State University,  \newline \indent 660041,
Krasnoyarsk, Russia.} \email{sadykov@lan.krasu.ru}
\thanks{Timur Sadykov was partially supported by the Russian
Foundation for Basic Research, grant 05-01-00517 and by the grant
MK-851.2006.1 of the President of Russian Federation.}

\begin{document}


\begin{abstract}
We consider algebraic functions~$z$ satisfying  equations of the
form
\begin{equation} \label{eq:1}
a_0 z^m + a_1z^{m_1} + a_2 z^{m_2} + \ldots + a_n z^{m_n} + a_{n+1} =0.
\end{equation}
Here $m > m_1 > \ldots > m_n>0,$ $m,m_i \in \N,$ and
$z=z(a_0,\ldots,a_{n+1})$ is a function of the complex variables
$a_0, \ldots, a_{n+1}.$ Solutions to such equations are
classically known to satisfy holonomic systems of linear partial
differential equations with polynomial coefficients. In this paper
we investigate one of such systems of differential equations which
was introduced by Mellin. We compute the holonomic rank of the
Mellin system as well as the dimension of the space of its
algebraic solutions. Moreover, we construct  explicit bases of
solutions in terms of the roots of~(\ref{eq:1}) and their
logarithms. We show that the monodromy of the Mellin system is
always reducible and give some factorization results in the
univariate case.

\end{abstract}


\maketitle

\section{Introduction}

\noindent
Consider the class of algebraic functions
which satisfy algebraic equations with symbolic coefficients,
i.e., equations of the form
\begin{equation}
a_0 z^m + a_1 z^{m_1} + a_2 z^{m_2} + \ldots + a_n z^{m_n} + a_{n+1} =0.
\label{generalAlgebraic}
\end{equation}
Here $m > m_1 > \ldots > m_n >0,$ $m,m_i \in \N,$
$z=z(a_0,\ldots,a_{n+1})$ is a function of the complex variables
$a_0, \ldots, a_{n+1}.$ Any solution to~(\ref{generalAlgebraic})
satisfies the $A$-hypergeometric system defined by Gelfand,
Kapranov and Zelevinsky~\cite{gkz89} (see also~\cite{Sturmfels,
Mayr}) associated to the matrix
$$ A =
\left(
\begin{array}{ccccc}
1 & 1& \dots &1 &1 \\
m & m_1 & \dots & m_n & 0
\end{array}
\right)
$$
and the homogeneity vector $(0, -1)$. In particular, a solution
to~(\ref{generalAlgebraic}) has a double homogeneity property and
can therefore be considered as a function of~$n$ variables. This
implies that we can arbitrarily prescribe the values of any two
nonzero coefficients in~(\ref{generalAlgebraic}) without losing
any essential information on the general solution to this equation.
It turns out to be convenient to divide~(\ref{generalAlgebraic})
by~$-a_{n+1}$ and then set $y=(-a_0 /a_{n+1})^{1/m}z,$
reducing~(\ref{generalAlgebraic}) to an equation of the form
\begin{equation}
y^{m} + x_1 y^{m_1} + \ldots + x_n y^{m_n} - 1 = 0.
\label{algebraic}
\end{equation}
A classical result of Mellin from~1921 (see~\cite{Mellin}) states
that solutions to such equations satisfy certain systems of
partial differential equations of hypergeometric type. More
precisely, if $y(x) = y(x_1,\ldots,x_n)$ is a solution
to~(\ref{algebraic}), then it satisfies the following system
of~$n$ partial differential equations:
$$
\prod_{k=0}^{m_j -1}
(m_1 \theta_1 + \ldots + m_n \theta_n + mk + 1)
\prod_{k=0}^{m_{j}^{'} -1}
(m_{1}^{'} \theta_1 + \ldots + m_{n}^{'} \theta_n + mk - 1) y(x) =
$$
\begin{equation}
(-1)^{m_j} m^m \frac{\partial^m y(x)}{\partial x_{j}^{m}},
\quad j = 1,\ldots, n,
\label{mellin}
\end{equation}
where $\theta_j = x_j\frac{\partial}{\partial x_j}$
and $m_{j}^{'} = m - m_{j}$.
In this article we will only consider algebraic equations of the
form~(\ref{algebraic}).
The system~(\ref{mellin}) will be referred to as
{\it the Mellin system of equations associated with~(\ref{algebraic}).}

The goal of this paper is to describe all the solutions of the Mellin system
of equations and to construct explicit bases of local
complex analytic solutions for any~$n$ and any value of the greatest
common divisor~$d$ of $m, m_1, \dots, m_n$. We describe
the space of algebraic solutions and we exhibit explicit non-algebraic
solutions in terms of the roots of~(\ref{algebraic}) and their logarithms.
These problems arose in the framework of the general theory of
hypergeometric functions and their relations with algebraic
equations, see~\cite{CDD},\cite{DMS},\cite{PassareTsikh},
\cite{Sturmfels}, and~\cite{walther}.

The special case of the trinomial equation $y^m + x y^p - 1 =0$
(which corresponds to ordinary Mellin equation) was thoroughly
studied in~\cite{KatoNoumi}.  The Mellin equation is always
reducible and  unless $p=m-1$, it is never the linear differential
equation with polynomial coefficients of  smallest order satisfied
by the roots of the trinomial equation. Moreover, in this case all
solutions to the Mellin equation are algebraic except in the case
when $d=1$ and $p < m-1$, where we exhibit an explicit
non-algebraic solution, as a special case of our main
Theorem~\ref{th:logs}.

We show that the monodromy of the Mellin equations is always
reducible. This could be deduced when $d=n=1$ from the general
results in~\cite{walther} since in this case the
homogeneity~$\beta$ has integer coordinates. In the last section,
we give some factorization results in the univariate case.

\medskip

\noindent{\bf Acknowledgements: }
The authors are grateful to Mart\'{\i}n Mereb for useful comments on
Theorem~\ref{th:logs} and to Michael Singer for interesting discussions
and for his help in the proof of Proposition~\ref{singer}.
T.~Sadykov is greatly indebted to his colleagues and staff at the
Max-Planck Institut f\"ur Mathematik in Bonn for creating a wonderful
working atmosphere.


\section{Convenient bases and generating solutions of the Mellin system}
\label{multidimsec}

Mellin calls the solution of~(\ref{algebraic}) which assumes the
value~$1$ for $x_1 = \ldots = x_n = 0$ {\it the principal
solution} of this equation. He shows (see~\cite{Mellin}) that the
power series expansion of the principal solution
to~(\ref{algebraic}) around the origin is given by
\begin{equation}
y_{pr}(x) = \sum_{\nu_1, \ldots, \nu_n \geq 0}
\frac{(-1)^{|\nu|}}{m^{|\nu|}} \,\, \frac {\prod_{\mu=1}^{|\nu|-1}
(m_1 \nu_1 + \ldots + m_n \nu_n - m\mu + 1)} {\nu_1! \ldots \nu_n!
} \,\, x_{1}^{\nu_1} \ldots x_{n}^{\nu_n}. \label{principalseries}
\end{equation}
Here $|\nu|=\nu_1+ \ldots + \nu_n$ and the empty product is
defined to be~$1.$

For a multi-index $I=(i_1,\ldots,i_n)$ and a complex vector $x =
(x_1,\ldots,x_n) \in \C^n$ we use the standard notation $x^I =
x_{1}^{i_1}\ldots x_{n}^{i_n}.$ By~$B$ we denote the set
$$
B \, = \,  \{
(i_1,\ldots,i_n) \in\Z^n : 0\leq i_j \leq m-1, {\rm\ for\ any\ }
j=1,\ldots, n \}.
$$

\begin{theorem} \label{th:convenient}
The dimension of the space of analytic solutions to~(\ref{mellin})
at a generic point, i.e., the holonomic rank of the Mellin system,
equals~$m^n$.  There exists a basis~$\{ f_{I}
\}_{I\in B}$ around the origin such that
\begin{equation} \label{Iform}
f_I(x) = x^I \tilde{f}_I (x_{1}^{m},\ldots,x_{n}^{m}),
\end{equation}
where~$\tilde{f}_I$ is analytic at~$0$ and $\tilde{f}_I(0) \neq 0$
for any~$I\in B.$ \label{structurethm}
\end{theorem}

We say that  a system of functions of the form~(\ref{Iform}) forms
a {\em convenient basis} in the solution space of the Mellin
system~(\ref{mellin}).

\medskip

The proof of this theorem can be deduced from Theorems~2.8 and~3.1
in~\cite{Sadykov} by adjusting the arguments in~\cite{Sadykov} to
the case of a nonconfluent hypergeometric system. For the benefit
of the reader we give now a proof which uses the results
of~\cite{Sadykov} as little as possible. We will reprove a part of
this theorem in the next section, as a consequence of the relation
with the $A$-hypergeometric system associated to the exponents
$m,m_1, \dots, m_n, 0$, and its translation to a Horn
system~\cite{DMS}.

\begin{proof}
For the sake of brevity we will use the notation $M = (m_1,
\ldots, m_n),$ $M^{'} = (m_{1}^{'}, \ldots, m_{n}^{'}),$
$\theta = (\theta_1, \ldots, \theta_n)$ and denote by~$\langle
\cdot,\cdot \rangle$ the scalar product. We multiply the $j$-th
equation in~(\ref{mellin}) with~$x_{j}^m.$ This will not affect
the space of holomorphic solutions to~(\ref{mellin}). Let us
denote by~$P_{j}(\theta)$ the polynomial $\prod_{k=0}^{m_j -1}
(\langle M,\theta \rangle + mk + 1) \prod_{k=0}^{m_{j}^{'}
-1} ( \langle M^{'}, \theta \rangle + mk - 1)$ and by~$G_j$
the differential operator
\begin{equation} \label{Gj}
G_j = x_{j}^{m} P_{j}(\theta) - (-1)^{m_j} m^m x_{j}^m
\frac{\partial^m}{\partial x_{j}^m}, \,\,\, j=1,\ldots, n.
\end{equation}
Let~$\mathcal{D}$ be the Weyl algebra in~$n$
variables~\cite{Bjork1}, and define $\mathcal{M}= \mathcal{D}/
\sum_{j=1}^n \mathcal{D}G_{j}$ to be the left $\mathcal{D}$-module
associated with the system of differential operators~(\ref{Gj}).
Let $R=\C[z_{1}, \ldots,z_{n}]$ and $R[x]=R[x_{1},\ldots,x_{n}]$
$=\C[x_{1},\ldots,x_{n},z_{1},\ldots,z_{n}].$ We make~$R[x]$ into a
left $\mathcal{D}$-module  by defining the action
of~$\partial_{j}$ on~$R[x]$ by
\begin{equation}
\partial_j = \frac{\partial}{\partial x_j} + z_{j}.
\label{derivation}
\end{equation}

Let us define the operators  $D_j : R[x] \rightarrow R[x]$ by
setting
\begin{equation}
D_j = z_j\frac{\partial}{\partial z_j} + x_{j}z_{j}, \,\,\,
j=1,\ldots,n. \label{ds}
\end{equation}
It was pointed out in~\cite{Adolphson} that the
operators~(\ref{ds}) form a commutative family of
$\mathcal{D}$-linear operators. Let~$D$ denote the vector
$(D_{1},\ldots,D_{n}).$ For any $j=1,\ldots,n$ we define the
operator $\nabla_{j}:R[x] \rightarrow R[x]$ by
$\nabla_{j}=z_{j}^{-1} D_{j}.$ This operator commutes with the
operators~$\partial_{k}$ since both~$D_{j}$ and the multiplication
by~$z_{j}^{-1}$ commute with~$\partial_{k}.$ Moreover, the
operator~$\nabla_{j}$  commutes with~$\nabla_{k}$ for all~$1\leq
j,k\leq n$ and with~$D_{k}$ for~$j\neq k.$ In the case~$j=k$ we
have $\nabla_{j} D_{j} = \nabla_{j}+D_{j}\nabla_{j}.$

Let us define a family of linear operators~$\{ W_j \}_{j=1}^n$
acting on the ring~$R[x]$ by
$$
W_j = P_j(D) \nabla_{j}^{m} - (-1)^{m_j} m^m \prod_{k=0}^{m-1}
(D_j - k), \,\,\, j = 1, \ldots, n.
$$
We claim that the family of operators~$\{ W_j \}_{j=1}^n$ is
commutative.
Indeed, using the fact that~$\{ D_j \}_{j=1}^n$ is a commutative
family of operators and the identity $\nabla_{j}^{m} D_j = (D_j +
m)\nabla_{j}^m$ we conclude, that
$$
W_j W_k = P_j(D) P_k(D_1,\ldots, D_j + m, \ldots, D_n) \,
\nabla_{j}^m \nabla_{k}^m =
$$
$$
\left( \prod_{\mu=0}^{m_j + m_k - 1} (\langle M, D \rangle +
m\mu + 1) \prod_{\mu=0}^{m_{j}^{'} + m_{k}^{'} - 1} ( \langle
M^{'}, D \rangle + m\mu - 1) \right) \nabla_{j}^m
\nabla_{k}^m.
$$
Computing the composition of the operators~$W_j$ and~$W_k$ in the
reversed order, we arrive at the same expression.

Using Lemma~2.2 and Theorem~2.4 in~\cite{Sadykov} we conclude that
the $\mathcal{D}$-module associated with the modified Mellin
system~(\ref{Gj}) is isomorphic to the quotient
$$
R[x] \Bigg/ \left( \sum_{j=1}^{n} W_j R[x]\right).
$$
By the commutativity of the operators~$W_j$ and Lemma~2.7
in~\cite{Sadykov}, the dimension of this quotient as a complex
vector space equals the product of the degrees of the operators in
the Mellin system, i.e.,~$m^n.$ This justifies the first claim of
the theorem.

Let us now prove the second statement of the theorem. Let $I\in B$
and consider a power series of the form
\begin{equation}
f_{I}(x) = \sum_{s\in m \N_{0}^{n} + I} \varphi(s) x^s,
\label{formalsolution}
\end{equation}
where $m \N_{0}^{n} + I = \{ (s_1, \ldots, s_n) \in\N_{0}^{n} :
s_k = mp_k + i_k, {\rm for\ some\ } p_k\in\N_{0}  \}.$ To compute
the action of the operator~$G_j$ on~$f_{I}(x)$ we notice that
\begin{equation}
x_{j}^{m} \frac{\partial^{m}}{\partial x_{j}^{m}} =
\prod_{k=0}^{m-1}(\theta_{j} - k) \label{useful}
\end{equation}
and that
$$
P_{j}(\theta) f_{I}(x) = \sum_{s\in m \N_{0}^{n} + I} P_j(s)
\varphi(s) x^s.
$$
Since $\prod_{k=0}^{m-1} (i_j - k)=0$ for any~$I\in B$ and any
$j=1,\ldots, n,$ it follows, that $G_{j} f_{I}(x) = 0$ if and only
if
\begin{equation}
\varphi(s) P_j(s + I) = (-1)^{m_j} m^m \varphi(s + m e_j)
\prod_{k=0}^{m-1} (s_j + m + i_j - k). \label{recurrencies}
\end{equation}
Here $e_j = (0,\ldots,1,\ldots,0)$ with~$1$ in the $j$-th place.
It is always possible to satisfy this relation by
choosing~$\varphi(s)$ in an appropriate way. Indeed, the product
$\prod_{k=0}^{m-1}(s_j + m + i_j - k)$ is never zero and hence the
recurrent relations~(\ref{recurrencies}) can always be solved.
Thus for any $I \in B$ there exists a formal series solution
to~(\ref{mellin}) of the form~(\ref{formalsolution}). Notice that
the resultant of the principal symbols of the operators in the
Mellin system~(\ref{mellin}) is different from zero at the origin.
By Proposition~2 in~\cite{PST}, any formal series solution to the
Mellin system of the form~(\ref{formalsolution}) converges in some
neighbourhood of the origin. The~$m^n$ series~$\{f_I(x)\}_{I\in
B}$ are linearly independent since their initial monomials are
different. Thus they form a convenient basis in the space of
analytic solutions to the Mellin system. This completes the proof
of Theorem~\ref{structurethm}.
\end{proof} 

\begin{definition}
\rm
We call an analytic solution~$y(x)$ to the Mellin system~(\ref{mellin})
defined near the origin {\it a generating solution} to this system if
for any~$I \in B$,
$$
\partial_I y(0) :=
\frac{\partial^{|I|}y(0)} {\partial x_{1}^{i_1} \ldots \partial
x_{n}^{i_n} } \neq 0.
$$
\label{gensoldef}
\end{definition}

As we will see in Theorem~\ref{basisthm}, a generating solution
gives rise to a basis in the solution space of the Mellin system.
On the other hand,  the sum of the elements of a convenient basis
is a generating solution to~(\ref{mellin}).

\begin{example}
\rm The Mellin equation associated with the quadratic equation
$y^2 + xy - 1 = 0$ has the form
\begin{equation}
(x^2 + 4)y'' + x y' - y = 0. \label{eq21}
\end{equation}
A convenient basis in the solution space of this second-order
differential equation is given by the functions $f_1 (x) = x ,$
$f_2 (x) = \sqrt{x^2 + 4}.$ Another basis is given by the roots
$y_i = \frac 1 2 ( - f_1  \pm f_2)$ of the quadratic equation $y^2+xy-1=0$,
which are generating solutions. \label{ex21}
\end{example}

We denote by~$\varepsilon$ the $m$-th
primitive root of unity $\varepsilon:=e^{2\pi i/m}$.
\begin{theorem}
Let~$y(x)$ be an analytic function near the origin which is a
solution to the Mellin system~(\ref{mellin}). Then,~$y(x)$ is a
generating solution if and only if the family of~$m^n$ functions
$y_{I}(x)= y(\varepsilon^{i_1}x_1, \ldots, \varepsilon^{i_n}x_n),$
$I= (i_1, \ldots, i_n) \in B$ is a basis in the solution space
of~(\ref{mellin}). \label{basisthm}
\end{theorem}
\begin{proof}
Assume that~$y(x)$ is a generating solution. The differential
operators $\theta_k=x_{k}\frac{\partial}{\partial x_k}$
and~$\frac{\partial^m}{\partial x_{k}^{m}}$ are invariant under
the change of variables $t_k = \varepsilon^{i_k} x_k$ for any
$k=1, \ldots, n.$ Indeed, $x_{k}\frac{\partial}{\partial x_k} =
\varepsilon^{-i_k}t_k \frac{\partial}{\partial t_k} \frac{\partial
t_k}{\partial x_k} = t_{k}\frac{\partial}{\partial t_k}.$ The
invariance of~$\frac{\partial^m}{\partial x_{k}^{m}}$ follows from
the identity~$\varepsilon^m = 1.$ Since the differential equations
in the system~(\ref{mellin}) only contain such operators, it
follows that the function $y_{I}(x) = y(\varepsilon^{i_1}x_1,
\ldots, \varepsilon^{i_n}x_n)$ solves~(\ref{mellin}) for any~$I\in
B.$

Let us show that the family of~$m^n$ functions~$\{ y_I(x) \}_{I\in
B}$ is linearly independent. This will prove the assertion of the
theorem. Let
$$
y(x) = \sum_{s\in\N_{0}^{n}} \varphi(s)x^s
$$
be the Taylor series expansion of~$y(x)$ around the origin. For a
given~$I\in B$ consider its subseries defined by
\begin{equation} \label{eq:subseries}
f_{I}(x) = \sum_{s\in m \N_{0}^{n} + I} \varphi(s) x^s,
\end{equation}
where $m \N_{0}^{n} + I = \{ (s_1, \ldots, s_n) \in\N_{0}^{n} :
s_k = mp_k + i_k, {\rm for\ some\ } p_k\in\N_{0}  \}.$ By
assumption~$y(x)$ is a generating solution to~(\ref{mellin}).
Hence the family of series~$\{ f_{I}(x) \}_{I\in B}$ is linearly
independent since their initial monomials are different. Being
subseries of a convergent power series, they all converge in some
neighbourhood of the origin. It follows directly from the
definition of~$f_{I}(x)$ that $y(x) = \sum_{I\in B} f_{I}(x).$

Suppose that the family of functions~$\{ y_{I}(x) \}_{I\in B}$
satisfies a linear relation
\begin{equation}
\sum_{I\in B} c_I y_I (x) \equiv 0 \label{relation}
\end{equation}
for some~$c_I \in \C.$ Let $J=(j_1, \ldots, j_n)\in B$ be a
multi-index. Using the definition of the function~$f_J$ and the
identity $\varepsilon^m = 1$ we conclude that
$$
f_{J}(\varepsilon^{i_1}x_1, \ldots, \varepsilon^{i_n}x_n) =
\sum_{s\in m\N_{0}^{n} + J} \varphi(s)
(\varepsilon^{i_1}x_1)^{s_1} \ldots (\varepsilon^{i_n}x_n)^{s_n} =
$$
$$
\sum_{s\in m\N_{0}^{n} + J} \varepsilon^{i_1 j_1 + \ldots + i_n
j_n} \varphi(s) x^s = \sum_{s\in m\N_{0}^{n} + J}
\varepsilon^{\langle I,J \rangle} \varphi(s) x^s=
\varepsilon^{\langle I,J \rangle} f_{J}(x),
$$
where $\langle I,J \rangle$ denotes the scalar product. Thus the
relation~(\ref{relation}) can be written in the form
$$
0 \equiv \sum_{I\in B} c_I y_I (x) = \sum_{I\in B} c_I
y(\varepsilon^{i_1}x_1, \ldots, \varepsilon^{i_n}x_n) =
$$
\begin{equation}
\sum_{I\in B} c_I \sum_{J\in B} f_{J}(\varepsilon^{i_1}x_1,
\ldots, \varepsilon^{i_n}x_n) = \sum_{J\in B} \left( \sum_{I\in B}
c_{I} \varepsilon^{\langle I, J \rangle} \right) f_{J}(x).
\label{relation1}
\end{equation}
The assumption that~$y(x)$ is a generating solution to the Mellin
system~(\ref{mellin}) implies that the family of functions~$\{
f_{J} \}_{J\in B}$ is linearly independent. It therefore follows
from~(\ref{relation1}) that
\begin{equation}
\sum_{I\in B} c_{I} \varepsilon^{\langle I, J \rangle} = 0, \,\,\,
{\rm for\ any\ } J\in B. \label{linsyst}
\end{equation}
The matrix of this system of linear homogeneous algebraic
equations is the $n$-th exterior power of the Vandermonde matrix
$$
\left(
\begin{array}{ccccc}
1 & 1                & 1             & \ldots & 1                   \\
1 & \varepsilon      & \varepsilon^2 & \ldots & \varepsilon^{m-1}   \\
1 & \varepsilon^2    & \varepsilon^4 & \ldots & \varepsilon^{2(m-1)}\\
\ldots & \ldots      & \ldots        & \ldots & \ldots              \\
1 & \varepsilon^{m-1}& \varepsilon^{2(m-1)}& \ldots & \varepsilon^{(m-1)^2}\\
\end{array}
\right)
$$
This Vandermonde matrix is nondegenerate since the numbers
$1,\varepsilon, \varepsilon^2, \ldots, \varepsilon^{m-1}$ are all
different. It is known that the eigenvalues of the exterior
product of two matrices are products of the eigenvalues of the
factors. Thus the exterior product of two nondegenerate matrices
is nondegenerate and hence so is the matrix of the
system~(\ref{linsyst}). This implies that~$c_{I}=0$ for any~$I\in
B$ and hence the functions~$\{ y_{I} (x) \}_{I\in B}$ are linearly
independent.

Assume now that~$y(x)$ is not generating and let~$J \in B$ such that
$\partial_J y(0) = 0$.  Then, $\partial_J y_I(0) = 0$ for all~$I \in B$. By
Theorem~\ref{th:convenient} there exists a solution~$f_J$ with
$\partial_J f_J(0) =1$.
Then, the family $\{y_I(x), \, I \in B\}$ cannot be a basis of solutions of the
Mellin system. This completes the proof.
\end{proof} 

Recall that we denote $ d = {\rm GCD}(m, m_1, m_2, \ldots, m_n) $.

\begin{corollary}\label{gcd}
If $d >1$ then $\{ y_{pr}(\varepsilon^{i_1}x_1, \ldots,
\varepsilon^{i_n}x_n) \}_{I\in B}$ is a basis in the solution
space to the Mellin system~(\ref{mellin}). Here~$y_{pr}(x)$ is the
principal solution~(\ref{principalseries}) to the Mellin system.
In particular, under the above assumption all solutions to the
Mellin system are algebraic functions. \label{prsolcor}
\end{corollary}

\begin{proof}
If $d >1$ then the term $m_1 \nu_1 + \ldots + m_n \nu_n - m\nu +
1$ in the numerator of the coefficient of~(\ref{principalseries})
is never zero. In particular, it is nonzero for $(\nu_1, \ldots,
\nu_n)\in B$ and hence~$y_{pr}(x)$ is a generating solution to the
Mellin system. The conclusion of the corollary follows now from
Theorem~\ref{basisthm}.
\end{proof} 

\begin{example}
\rm
Consider the algebraic equation
\begin{equation}
y^6 + x_1 y^4 + x_2 y^2 - 1 = 0.
\label{eq642}
\end{equation}
A generating solution for the corresponding Mellin system is given
by the function
$$
y(x_1,x_2) =
\sqrt{ -2 x_1  + u(x_1,x_2) - \frac{12 x_2 - 4 x_{1}^2}{u(x_1,x_2)}},
$$
where
$
u(x_1,x_2) = \sqrt[3]{108 + 36 x_1 x_2 - 8 x_{1}^3 +
12\sqrt{81 + 54 x_1 x_2 - 12 x_{1}^3 + 12 x_{2}^3 - 3 x_{1}^2 x_{2}^2}}.$
A simple computation in Maple shows that
all of the initial exponents for~(\ref{eq642}) are present
in the series expansion of~$y(x_1,x_2)$, as predicted by Corollary~\ref{gcd}.
In contrast with Example~\ref{ex321} below, in this case it is possible
to give a generating solution to the Mellin system associated
with~(\ref{eq642}) in terms of the solutions to~(\ref{eq642}) itself.
\label{ex642}
\end{example}

\begin{remark}
\rm From the point of view of the original algebraic
equation~(\ref{algebraic}) it is, of course, very unnatural to
make the assumption that $d >1.$ Indeed, the change of  unknown
$u(x) = y(x)^{d}$ reduces~(\ref{algebraic}) to an equation with
relatively prime exponents. However, it turns out that the Mellin
system associated with this transformed equation differs
essentially from the Mellin system associated with the original
algebraic equation. This is illustrated by Examples~\ref{62}
and~\ref{ex31}, and by Examples~\ref{ex642} and~\ref{ex321},
since in  case $d>1$ all solutions are
algebraic by Corollary~\ref{gcd}, while as we will see, in  case
$d=1, n=1, m_1 < m-1$ or
$d=1, n >1$, there always exists a non-algebraic solution.
Thus from the point of view of
the Mellin system the case when $d>1$ cannot be trivially reduced
to the case of relatively prime exponents.
\label{gcdrem}
\end{remark}


\section{$A$-hypergeometric, Horn hypergeometric and Mellin systems}
\label{Asec}

Given the equation~(\ref{generalAlgebraic}), its roots (thought of as
functions of $a = a_0, \dots, a_{n+1})$) satisfy the
$A$-hypergeometric system with parameter~$\beta := (0, -1)$
(cf~\cite{Sturmfels}), where
\begin{equation} \label{matrixA}
A := \, \left(
\begin{array}{ccccc}
1 & 1& \dots &1&1 \\
m & m_1 & \dots & m_n & 0
\end{array}
\right).
\end{equation}
Namely, it is the left ideal in the Weyl algebra
$\C[a_0,\dots,a_{n+1},\pd 0,\dots,\pd {n+1}]$ generated by
\begin{eqnarray}
\label{eq:higherorder} \,\,\,\, \hbox{the toric operators} \quad
\partial^u - \partial^v \quad \hbox{for}\quad
 u,v\in \N^{n+2}   \quad \hbox{with}\quad A\cdot u=A\cdot v, \\
\label{eq:firstorder} \,\, \hbox{and the Euler operators} \quad
\sum_{j=0}^{n+1} a_j \partial_j \quad \hbox{and} \quad m a_0 +
\sum_{j=1}^{n} m_j a_j \partial_j +1. \phantom{wow}
\end{eqnarray}

Note that the corresponding $A$-hypergeometric system is the same
as the hypergeometric system associated with the matrix
\begin{equation} \label{matrixA'}
A' := \, \left(
\begin{array}{ccccc}
1 & 1& \dots &1&1 \\
m/d & m_1/d & \dots & m_n/d & 0
\end{array}
\right)
\end{equation}
and homogeneity $ \beta'=(0, -1/d).$

Consider now the following matrix $\mathcal{B} \in \Z^{(n+1)
\times m}$:
\begin{equation} \label{matrixB}
\mathcal{B} := \, \left(
\begin{array}{cccc}
-m_1& -m_2& \dots &-m_n \\
m & 0 & \dots & 0 \\
0 & m & \dots & 0 \\
\vdots& \vdots&\vdots & \vdots \\
0 & 0 & \dots & m\\
-m'_1 & - m'_2 &\dots & - m'_n
\end{array}
\right),
\end{equation}
where as above $m'_i = m - m_i, \, i=1, \dots, n$. Let $c =(-1/m,
0 \dots, 0, 1/m)$ and denote by ${\rm Horn}\,_{\mathcal B}(c) =
\langle H_1, \dots, H_n \rangle$ the corresponding Horn system, as
in Definition~2.1 in~\cite{DMS}, which we now recall.
For $j=1, \dots, n,$ the Horn
operator~$H_j$  in the variables $w = (w_1, \dots, w_n)$ equals
\begin{equation} \label{eq:Hornop}
 \prod_{k=0}^{m-1} (m \theta_j -k) - w_j
\prod_{k=0}^{m_j -1} (-m_1 \theta_1 - \ldots - m_n \theta_n -
\frac 1 m -k) \prod_{k=0}^{m_{j}^{'} -1} (-m_{1}^{'} \theta_1 -
\ldots - m_{n}^{'} \theta_n + \frac 1 m - k).
\end{equation}

We deduce from Lemma 5.1 and Corollary 5.2 in~\cite{DMS} the
following result.

\begin{lemma} \label{Ahorn}
Let $\psi(w_1, \dots, w_n)$ be a function of~$n$ complex variables
and let~$\varphi(a)$ be the following function of~$n+2$ variables
(defined in suitable simply connected open sets):
\begin{equation} \label{eq:phipsi}
\varphi(a) \, = \, \left( \frac {a_0}{a_{n+1}} \right)^{- \frac  1
m} \psi \left( \frac{a_1^m}{a_0^{m_1} a_{n+1}^{m'_1}}, \dots,
\frac{a_n^m}{a_0^{m_n} a_{n+1}^{m'_n}}\right).
\end{equation}
Then,~$\varphi$ is a solution to the $A$-hypergeometric system
with parameter~$\beta$ if and only if~$\psi$ is a solution to
${\rm Horn}\,_{\mathcal B}(c)$.
\end{lemma}

Setting $a_0 =1, a_{n+1} = -1$, we deduce the relation
$$ \varphi(1, a_1, \dots, a_n, -1) = (-1)^{\frac 1 m} \,
\psi ((-1)^{m'_1} a_1^m, \dots,
(-1)^{m'_n} a_n^m).$$ We thus make the change of variables $w_j =
(-1)^{m'_j} x_j^m$ and we get that given a holomorphic function
$\psi(w_1, \dots, w_n)$ as above, the function
\begin{equation}\label{eq:mupsi}
\mu (x_1, \dots, x_n):= \, \psi(w(x)) \, = \, \psi((-1)^{m'_1}
x_1^m, \dots, (-1)^{m'_n} x_n^m)
\end{equation}
verifies
\begin{equation}\label{eq:muvarphi}
(-1)^{\frac 1 m} \mu(a_1, \dots, a_n) = \varphi(1, a_1, \dots,
a_n, -1).
\end{equation}

We now translate the Horn operators in the~$w$ variables to
the~$x$ variables. Since $w_j \frac \partial {\partial w_j} =
\frac 1 m x_j \frac \partial {\partial x_j}$, we get for all $j=1,
\dots, m$,
$$
H'_j := \, \prod_{k=0}^{m-1} (\theta_j -k) - (-1)^{m'_j} x_j^m
\prod_{k=0}^{m_j -1} \left( -\frac{m_1}{m} \theta_1 - \ldots -
\frac{m_n}{m} \theta_n - \frac 1 m -k \right)\circ
$$
\begin{equation}\label{eq:horntomellin1}
\prod_{k=0}^{m_{j}^{'} -1} \left( -\frac{m_{1}^{'}}{m} \theta_1 -
\ldots - \frac{m_{n}^{'}}{m} \theta_n + \frac 1 m - k \right),
\end{equation}
where now $\theta_1, \dots, \theta_n$ refer to derivatives with
respect to the variables $(x_1, \dots, x_n)$, i.e., $\theta_i =
x_i \frac{\partial}{\partial x_i}, \, i=1, \dots, n$.

Note that $\prod_{k=0}^{m-1} (\theta_j -k) = x_j^m \left(\frac
\partial {\partial {x_j}}\right)^m$. Multiplying~$H'_j$ by
$(-1)^{m+1} m^m$ and dividing by~$x_j^m$ we get the Mellin
operators~(\ref{mellin}).

We immediately deduce:

\begin{lemma}\label{lemma:mupsi}
A holomorphic function $\mu(x_1, \dots, x_n)$ is annihilated by
the Mellin operators~(\ref{mellin}) if and only if the function
$\psi(w_1,\dots, w_n)$ verifies the Horn system~(\ref{eq:Hornop}),
where~$\mu$ and~$\psi$ are related by~(\ref{eq:mupsi}) (in
suitable open sets).
\end{lemma}

We now take advantage of these translations to provide a new proof
of the rank of the Mellin system stated in
Theorem~\ref{structurethm}. Note that the greatest common
divisor~$g$ of the maximal minors of~$B$ equals $g = m^{n-1} d$. The
normalized volume of~$A$ equals~$m/d$. Moreover, by  the results
in~\cite{fs}, the lattice ideal~$I_B$ is a complete intersection.
 We deduce from~\cite[Theorem~2.5]{DMS}
that the rank of ${\rm Horn}\,_{\mathcal B}(c)$ equals
$m^{n-1} d \times m/d = m^n$. Using the translation of solutions in
Lemma~\ref{lemma:mupsi}, it is easy to check that the holonomic
rank of the Mellin system is also~$m^n$.

\begin{theorem} \label{th:phimu}
In suitable open sets, given an $A$-hypergeometric
 function~$\varphi(a)$ with parameter~$\beta= (0,-1)$,
the function of~$n$ variables defined by specialization
\begin{equation}\label{eq:esp}
\mu(x_1, \dots, x_n) = \varphi( 1, x_1, \dots, x_n, -1),
\end{equation}
is a solution to the Mellin system.
\end{theorem}

\begin{proof}
Given a choice of $m$-roots of the coordinates, define
$$
\psi (w_1, \dots, w_n):= \mu(\varepsilon^{-m'_1} w_1^{\frac 1 m},
\dots, \varepsilon^{-m'_n} w_n^{\frac 1 m}),
$$
where~$\varepsilon'$
is a primitive $m$-root of~$-1$. Note that $w_i = (-1)^{m'_i} \left(
{\varepsilon'}^{-m'_i} w_i^{\frac 1 m}\right)^m, \, i=1,\dots, n$. By
Lemma~\ref{lemma:mupsi}, it is enough to verify that~$\psi$ is a
solution to the Horn system. Now, by Lemma~\ref{Ahorn}, this
happens if and only if the function~$\varphi'$ defined by:
$$\varphi'(a) :=  \, \left( \frac {a_0}{a_{n+1}}
\right)^{- \frac 1 m} \psi \left( \frac{a_1^m}{a_0^{m_1}
a_{n+1}^{m'_1}}, \dots, \frac{a_n^m}{a_0^{m_n}
a_{n+1}^{m'_n}}\right),$$ satisfies the $A$-hypergeometric system
with parameter~$\beta$. But calling~$\lambda$ an $m$-th root of $-
a_0/a_{n+1}$ we deduce from the definitions of the functions and
the homogeneity $(0,-1)$ satisfied by~$\varphi$ that
$$
\begin{array}{ll}
\varepsilon' \, \varphi'(a) \, &=
\, \varepsilon' \, \left( \frac {a_0}{a_{n+1}} \right)^{- \frac 1 m} \,
\psi
(a_{1}^m a_{n+1}^{m_1 - m} a_{0}^{-m_1}, \dots,
a_{n}^m a_{n+1}^{m_n - m}  a_{0}^{-m_n})  \\
&= \, \lambda^{-1}  \mu ((- \frac {a_1}{a_{n+1}}) \lambda^{-m_1} ,
\dots,
 (- \frac {a_n}{a_{n+1}}) \lambda^{-m_n}) \\
&= \, \lambda^{-1} \varphi(1, - \frac {a_1}{a_{n+1}}\lambda^{-m_1},
\dots, - \frac {a_n}{a_{n+1}}\lambda^{-m_n}, -1)= \varphi(a),
\end{array}
$$
and we are done.
\end{proof}


\section{Solutions in terms of roots} \label{rootsec}

Mellin not only observed in~\cite{Mellin} that the roots $y_1(x),
\dots, y_m(x)$ of the algebraic equation~(\ref{algebraic})
$$
y^{m} + x_1 y^{m_1} + \ldots + x_n y^{m_n} - 1 = 0
$$
satisfy the Mellin system~(\ref{mellin}), but he also made the
following easy observation. Given the principal solution
$y_{pr}(x)$  in~(\ref{principalseries}), all the solutions have
the form
$$ \eta \, y_{pr} (\eta^{m_1} x_1, \dots, \eta^{m_n} x_n),$$
where~$\eta$ runs through the $m$-roots of~$1$.

It is also clear that for any choice of $I=(i_1,\ldots,i_n) \in
\N^n$, the function $y_{I}(x)= y_{pr}(\varepsilon^{i_1}x_1,
\ldots, \varepsilon^{i_n}x_n)$ is a root of the algebraic equation
$$
(I) \phantom{---------} y^{m} +  \e^{i_1} x_1 y^{m_1} + \ldots +
\e^{i_n} x_n y^{m_n} - 1 = 0, \phantom{--------------}
$$
where, as in Section~\ref{multidimsec}, we denote $\e =  e^{2\pi
i/m}$. Moreover, as we have already remarked, all these functions
lie in the solution space of~(\ref{mellin}). So,  the Mellin
system has not only the roots of the algebraic
equation~(\ref{algebraic}) as solutions, but also the roots of the
associated equations~$(I)$. It is clear that we only get~$m^n$
different equations this way, and we could parametrize them
taking~$I\in B$. We denote by~$Y$ the $\C$-vector space generated
by all the roots of all the equations of the form~$(I)$.   The
roots of the original equation~(\ref{algebraic})  correspond
precisely to the functions $\e^j y_{(j m_1, \dots, j m_n)}, \,
j=0, \dots, m-1$. Note that these~$m$ functions are all distinct,
but of course they are not linearly independent if $d > 1$ or if
$d=1$ and $m_1 < m-1$. In case $d>1$ we have seen in
Corollary~\ref{gcd} that~$Y$ coincides with the solutions space of
the Mellin system.

\begin{example}
The polyquadratic equation.
\rm
Consider the equation
\begin{equation}
y^{2k} + xy^k - 1 = 0.
\label{eq2kk}
\end{equation}
Its~$2k$ roots are $\varepsilon^j \sqrt[k]{(-x \pm \sqrt{x^2 + 4})/2},$
$j=0,\ldots,k-1,$ where~$\varepsilon$ is the primitive $k$-th root
of unity.
Observe that these roots span a vector space of dimension~$2.$
A basis in the $2k$-dimensional solution space of the Mellin equation
associated with~(\ref{eq2kk}) is given by the functions
$\sqrt[k]{-x \pm \sqrt{x^2 + 4\varepsilon^j}},$ $j = 0,\ldots, k-1.$
\label{ex2kk}
\end{example}

For the rest of this section we will then assume that $d=1$.
Our task is to explain and generalize the following example.

\begin{example}
The general cubic.
\rm
Consider the algebraic equation
\begin{equation}
y^3 + x_1 y^2 + x_2 y - 1 = 0.
\label{eq321}
\end{equation}
The corresponding Mellin system is given by
\begin{equation}
\begin{array}{rcl}
 27 \theta_1 (\theta_1 - 1)(\theta_1 - 2)y &=&
x_{1}^{3} (2\theta_1 + \theta_2 + 1)(2\theta_1 + \theta_2 + 4)
(\theta_1 + 2\theta_2 - 1)y, \\
-27 \theta_2 (\theta_2 - 1)(\theta_2 - 2)y &=&
x_{2}^{3} (2\theta_1 + \theta_2 + 1)(\theta_1 + 2\theta_2 - 1)
(\theta_1 + 2\theta_2 + 2)y.
\end{array}
\label{mellin321}
\end{equation}
The three solutions to~(\ref{eq321}) are linear combinations of
the following functions:
$$
y_1(x_1,x_2) = x_1, \quad
y_2(x_1,x_2) = \sqrt[3]{u(x_1,x_2)}, \quad
y_3(x_1,x_2) = \frac{3 x_2 - x_{1}^{2}}{\sqrt[3]{u(x_1,x_2)}},
$$
where
$u(x_1,x_2) = 108 + 36 x_1 x_2 - 8 x_{1}^{3} + 12\sqrt{81 + 54 x_1 x_2
- 12 x_{1}^{3} + 12 x_{2}^{3} - 3 x_{1}^{2} x_{2}^{2}}.$
The expansions of~$y_2$ and~$y_3$ into power series around the origin
are given by
$$
y_2 =
6 + \frac{2}{3} x_1 x_2 - \frac{4}{27} x_{1}^{3} +
\frac{2}{27} x_{2}^{3} - \frac{4}{27} x_{1}^{2} x_{2}^{2} + \ldots, \quad
y_3 =
\frac{1}{2} x_2 - \frac{1}{6} x_{1}^{2} - \frac{1}{18} x_{1} x_{2}^{2}
+ \ldots
$$
In this example, the set of initial exponents is given by
$B=\{ (i,j)\in\Z^2 : 0 \leq i,j \leq 2 \}.$ The expansions of~$y_1,$
$y_2,$ and~$y_3$ into power series contain~$7$ of the~$9$ elements
in~$B.$  This could be also checked differentiating ~(\ref{eq321}) implicitly
and evaluating at~$0$.
The~$2$ missing exponents are~$(2,1)$ and~$(0,2).$
Thus the solution space of~(\ref{mellin321}) contains a two-dimensional
subspace which is not spanned by solutions to~(\ref{eq321}).
\label{ex321}
\end{example}

Remark  that the roots of the equations~$(I)$  above with any
index set of the form $I_j = (jm_1, \dots, jm_n)$ for any $j=1,
\dots, m-1$, are just constant multiples of the roots of the
original equation~(\ref{algebraic}) (corresponding to $j=0$). More
generally, there are $m^{n-1}$ subsets $G^{(1)}, \dots,
G^{(m^{n-1})}$ of  algebraic equations of the form~$(I),$ each of
them consisting of~$m$ equations, such that for any pair of
equations~$(I), (J)$ in the same subset~$G^{(k)}$, the
corresponding solutions only differ by multiplicative constants.
Denote by $\Gamma  = \{ I^{(1)}, \dots, I^{(m^{n-1})}\},$ where
each $I^{(k)} \in B$ is  a choice of an index set in the subset
$G^{(k)}$ of algebraic equations, for all $k= 1, \dots, m^{n-1}$.
The analytic  roots of the algebraic equation $(I^{(k)})$ defined
in a small neighborhood of the origin will be denoted by
$y^{(k)}_1, \dots, y^{(k)}_m$.

Given a complex vector $c =(c^{(1)}, \dots, c^{(m^{n-1})}),$
denote by~$\chi_c$ the function defined by
\begin{equation}\label{eq:chi}
\chi_c \, = \, \sum_{k=1}^{m^{n-1}} \,
c^{(k)}  (y^{(k)}_1 \log(y^{(k)}_1) +\dots + y^{(k)}_m \log(y^{(k)}_m)),
\end{equation}
where~$\log$ denotes a holomorphic branch of the logarithm defined
in the union of the images of the roots.

When the exponents $m,m_1, \ldots, m_n$ are relatively prime,
is it shown in~\cite[Th. 2.4]{CDD} (see also~\cite[Cor. 3.6]{KatoNoumi})
that  the dimension of the vector space generated by the solutions
to the algebraic equation~(\ref{algebraic}) is equal to~$m$ if
$m_1 = m-1$ and it is equal to~$m-1$ if $m_1 < m-1.$
Thus, when $d=1$ the only possible linear relation
between the solutions to the algebraic equation~(\ref{algebraic})
is the obvious one: the roots of~(\ref{algebraic}) sum up to zero
if the coefficient by~$y^{m-1}$ in~(\ref{algebraic}) is zero.

Our main result is the following.

\begin{theorem} \label{th:logs}
The space~$Y$ coincides with the space of algebraic solutions to
the Mellin system of equations around the origin. Its dimension
equals $\dim Y = m^n - m^{n-1} +1$ in case $m_1=m-1$ and $\dim Y =
m^n - m^{n-1}$ otherwise. The space of relations
$$R := \{ c \in \C^{m^{n-1}} \,
 /  \,
\sum_{k=1}^{m^{n-1}} \, c^{(k)}  (y^{(k)}_1  +\dots + y^{(k)}_m)  = 0 \}$$
has dimension $\dim R = m^{n-1} -1$ in case $m_1 = m-1$ and $\dim R = m^{n-1}$
otherwise. Moreover, for any nonzero $c \in R$, the function~$\chi_c$
defined by~(\ref{eq:chi})
is a non-algebraic solution to the Mellin system, the space
$$
S := \{ \chi_c \, / \, c \in R\}
$$
has dimension $\dim S = \dim R$ and the sum $S \oplus Y$ is direct and
equals the full solution space of the Mellin system.
\end{theorem}

When $n=1$, $S$ is nonzero only when $m_1 < m-1$ and in this case
it is spanned by the non-algebraic function $\chi = y_1 \log(y_1)
+ \dots + y_m \log(y_m)$, where $y_1(x), \dots, y_m(x)$ denote the
roots of $y^m + x y^{m_1} -1 =0$, while the space of algebraic
solutions~$Y$ has dimension $m-1$. This is the content of
Theorem~3.5 in~\cite{CDD} for $A$-hypergeometric systems
associated to monomial curves, which can be directly translated to
the Mellin equation via our computations in Section~\ref{Asec}.
So, Theorem~\ref{th:logs} is inspired by and generalizes this
result. When translating the Mellin system to  the homogeneous
setting as in Section~\ref{Asec}, the ocurrence of logarithmic
solutions is explained by the results in~\cite{dmm}.

\begin{example} \label{ex321cont}
{\bf (Example~\ref{ex321} continued)} \rm In this case $\dim R =
\dim S = 3 -1 =2$, as we have seen. There are~$3$ subsets of~$3$
equations each and we could choose the following representatives
to form our set~$\Gamma$:
$$
\begin{array} {lllr}
(I_1) \phantom{-----------} y^3 + x_1 y^2 + x_2 y - 1 &=& 0,   \phantom{-------------} \\
(I_2) \phantom{-----------} y^3 + x_1 y^2 +  \e x_2 y - 1 &=& 0,  \phantom{-------------} \\
(I_3) \phantom{-----------} y^3 + x_1 y^2 + \e^2 x_2 y - 1 &=& 0, \phantom{-------------}
\end{array}
$$
where~$\e$ is a primitive cubic root of~$1,$  and $I_1 = (0,0),$
$I_2 = (0,1),$ $I_3 =  (0,2).$

A basis of relations in~$R$ is given for instance by the  two
relations $(1,-1,0),$ $(1, 0, -1)$, since the sum of the roots of
each of the~$3$ equations above equals the same value~$x_1$. On
the other side, the space~$Y$ spanned by the functions $y^{(k)}_j,
\, j,k = 1,2,3$ has dimension~$7$. So a basis of the solutions
near the origin is given by choosing seven of these roots together
with the two non-algebraic functions $\chi_{(1,-1,0)},$
$\chi_{(1,0,-1)}$.
\end{example}

\noindent {\it Proof of Theorem~\ref{th:logs}.} We begin by
computing the dimension of~$Y$. Call~$\varphi_\nu$, for $\nu \in
\N_0^n$, the coefficients in the expansion~(\ref{principalseries})
of the principal solution $y = y_{pr}$ of the algebraic
equation~(\ref{algebraic}):
\begin{equation}
\varphi_\nu \, = \, \frac{(-1)^{|\nu|}}{m^{|\nu|}} \,\, \frac {\prod_{\mu=1}^{|\nu|-1}
(m_1 \nu_1 + \ldots + m_n \nu_n - m\mu + 1)} {\nu_1! \ldots \nu_n!}.
\end{equation}
Note that when $d>1,$ the solution~$y_{pr}$ is generating. This is
proven in Theorem~\ref{basisthm}, whose proof we now generalize
with similar notations.

Let $f_I, I \in B,$ be defined as in~(\ref{eq:subseries}). Observe
that when $\varphi_I =0$, then $\varphi_{I'}=0$ for all $ I' \in
I + m \N_{0}^{n}$. It follows that when $\varphi_I=0$ then
$f_I=0$. Let $B' = \{ I \in B \, : \, \varphi_I \not=0\}.$ We
claim that $ \dim Y = \# B'$.   Clearly, $ y = \sum_{I \in B'}
f_I$ and the family $\{f_I\, : \, I \in B'\}$ is linearly
independent. Suppose that the family of functions $\{y_I\}_{I \in
B}$ satisfies a linear relation
\begin{equation} \label{eq:equiv0}
\sum_{I\in B} c_I y_I (x) \equiv 0,
\end{equation}
for some~$c_I \in \C.$  Then,
$$
\sum_{J\in B'} \left( \sum_{I\in B}
c_{I} \varepsilon^{\langle I, J \rangle} \right) f_{J}(x) = 0,
$$
which implies that
\begin{equation} \label{eq:B'}
\sum_{I\in B} c_{I} \varepsilon^{\langle I, J \rangle} = 0, \,\,\,
{\rm for\ any\ } J\in B'.
\end{equation}
In fact,~$(c_I)_{I \in B}$ satisfy~(\ref{eq:equiv0}) if and only
if these last equalities hold. Since the full  matrix
$(\varepsilon^{\langle I, J \rangle} )_{I,J \in B}$ is of maximal
rank, it follows that the system of equations~(\ref{eq:B'}) has
rank~$\# B'$ and so the dimension of the space of linear
relations~$R'$ among the algebraic functions $\{ y_I, \, I \in
B\}$ is equal to $\#B - \# B'$. Dualizing, we conclude that their
linear span~$Y$ has dimension~$\# B'$, as wanted.

We now prove that the cardinality of~$B'$ equals $m^n - m^{n-1}
+1$ in case $m_1=m-1$ and $\dim Y = m^n - m^{n-1}$ otherwise. Call
$B'' = B \setminus B'$ its complement in~$B$. It is clear
that~$\nu \in B''$ if and only if there exists a positive
integer~$\mu$ strictly smaller than~$|\nu|$ such that
$$
m_1 \nu_1 + \dots + m_n \nu_n - m \mu + 1 =0,
$$
which we simply write as $ \langle M, \nu \rangle - m \mu + 1 =0,$
where $M = (m_1, \dots, m_n).$ Note that when this equality holds,
it follows that the integer~$\mu$ satisfies $1 \leq \mu \leq
|\nu|$, since $0 < m_n < \dots < m_1 < m$  and all~$\nu_i$ are
non-negative. Moreover, it is easy to check that the equality $\mu
= |\nu|$ can only hold when $m_1=m-1$ and $ \mu = \nu_1 =1$ and
$\nu_j=0$ for any $j >1$. So it is enough to prove that the
cardinality of the set
$$
\{ \nu \in B \, : \,  \langle M , \nu \rangle  = -1  \mod(m)\}
$$
equals $m^{n-1}$.  We claim that for any  $r = 0, \dots, m-1$
there exist precisely $m^{n-1} = m^n/m$ solutions $\nu\in B$ to
the modular equation $ \langle M , \nu \rangle = r  \mod(m)$.
Since we are assuming that $d=1$, there exists an integer
$n$-tuple $X=(x_1, \dots, x_n)$ satisfying $ \langle M, X \rangle
= 1 \mod(m)$ and so the equation $\langle M, \nu \rangle = r
\mod(m)$ has a solution for any~$r$. Then it is easy to see that
for any~$\nu$, $\{\langle M, \nu \rangle \mod(m) , \langle M, \nu
+ X \rangle \mod(m),  \dots, \langle M, \nu + (m-1)X\rangle
\mod(m)\}$ takes all the values $r =0, \dots, m-1$ precisely once.
Moreover, two such subsets are disjoint or equal, from which the
result follows since the cardinality of~$B$ is~$m^n$.

It follows that the space~$R'$ of linear relations among the
algebraic functions $\{ y_I, \, I \in B\}$ has dimension
$m^{n-1}-1$ in case $m_1=m-1$ and $m^{n-1}$ otherwise. The
relations~$R$ define in principle a subspace of~$R'$ corresponding
to those vectors~$c$ which have a special pattern of repeated
coordinates, but it is easy to see that it has the same dimension.
Indeed it is very easy to find a linear relation between the sum
of the roots in each equation (which equals  minus the coefficient
of~$y^{m-1}$, which is a constant multiple of~$x_1$ if $m_1 =m-1$
and~$0$ otherwise), with the sum of the roots of the original
equation~(\ref{algebraic}). This gives $m^{n-1} -1$ independent
linear relations in~$R$ and we get moveover the relation $c= (1,
0, \dots, 0)$ in case $m_1 < m-1.$ We can explicitly list a basis of the
space~$R$: the vectors  $\{e_1-e_2,...,e_1-e_{m^{n-1}}\}$ when $m_1 = m-1$
and  the vectors $\{e_1,...,e_{m^{n-1}}\}$ otherwise.

Then,~$R=R'$ and we recover in
particular the quoted result in~\cite[Th. 2.4]{CDD} about the only
possible linear relation among the roots of a single sparse
algebraic equation with generic coefficients.

In order to see that for any $c \in R$ the function~$\chi_c$ is a
solution to the Mellin system, we will need the translations
described in Section~\ref{Asec} between the Mellin system and the
corresponding Horn and $A$-hypergeometric systems. This will allow
us to use the results in~\cite{CDD} based on the study of the
relation between lattice hypergeometric ideals and
$A$-hypergeometric ideals explained in~\cite{DMS}. We sketch the
main ingredients of the proof and refer the reader to~\cite{CDD,
DMS} for the details. Since the lattice generated by the columns
of the matrix $\mathcal{B} \in \Z^{(n+1) \times m}$
in~(\ref{matrixB}) equals $m^{n-1}$ and the associated lattice
ideal is a complete intersection, the Horn
system~(\ref{eq:Hornop}) gets translated by homogenization (as in
Lemma~\ref{Ahorn}) to the hypergeometric system ${\mathcal
H}_{\mathcal B}$ corresponding to this lattice. We get~$m^{n-1}$
hypergeometric systems ${\mathcal H}_i, i = 0, \dots, m^{n-1}-1$
which are torus translates of the $A$-hypergeometric system, all
with the same homogeneity $\beta= (0,-1)$ and such that all
solutions of any~${\mathcal H}_i$ satisfy the system ${\mathcal
H}_{\mathcal B}$. The different functions $y_I \in Y$ which
satisfy the Mellin system, and thus the Horn system after the
change of variables in Lemma~\ref{lemma:mupsi}, satisfy (after
$A$-homogenization as in Lemma~\ref{Ahorn}) different systems in
the family ${\mathcal H}_i, i = 0, \dots, m^{n-1}-1$ according
to~$(I).$ Each system is satisfied precisely by~$m$ of these
homogenized functions, which we denote by~$z_I(a)$. The
functions~$z_I(a)$ also satisfy the corresponding algebraic
equation:
\begin{equation} \label{eq:I}
a_m z^m + \e^{i_1} a_1z^{m_1} + \e^{i_2} a_2 z^{m_2} + \ldots +
\e^{i_n} a_n z^{m_n} + a_{n+1} =0.
\end{equation}
In particular,
when $I=(i_1, \dots, i_n)$ with $i_k = j m_k \mod(m)$ for all $k=1,\dots, n$
and some fixed $j =0, \dots, m-1$,
the function~$z_I$  satisfies the differential $A$-hypergeometric
system recalled in Section~\ref{Asec} and the algebraic equation~(\ref{generalAlgebraic}).
Recall that $y_I(x) = z_I(1, x_1, \dots, x_n,-1)$ as in Theorem~\ref{th:phimu}.
Changing a little bit the notation of the homogeneous roots
as in~(\ref{eq:chi}), in order to prove that  the function~$\chi_c$ is a
solution of the Mellin system, it is enough to prove that the homogenized function
\begin{equation} \label{eq:overlinechi}
\overline{\chi}_c:= \, \, \sum_{k=1}^{m^{n-1}} \,
c^{(k)}  (z^{(k)}_1 \log(z^{(k)}_1) +\dots + z^{(k)}_m \log(z^{(k)}_m))
\end{equation}
satisfies the lattice hypergeometric system ${\mathcal
H}_{\mathcal B}.$ Then, it is possible to prove that  each of the
functions 
\begin{equation} \label{eq:alpha}
\alpha^{(k)}:= z^{(k)}_1 \log(z^{(k)}_1) +\dots + z^{(k)}_m
\log(z^{(k)}_m)
\end{equation}
is annihilated by the toric operators in the
corresponding torus translate of the $A$-hypergeometric system as
in the proof of Proposition~3.2 in~\cite{CDD}. Moreover, the fact
that~$c$ lies in~$R$ allows one to prove that the whole sum
$\overline{\chi}_c$ satisfies the homogeneity equations.

The theorem will be proved if we show the following statement:
If a linear combination $\chi_\lambda:=\sum_{i=1}^s \lambda_i \chi_{c[i]}$ is
algebraic, where $\lambda_i \in \C$ and $c[1], \dots, c[s] \in R$,
then $\sum_{i=1}^s \lambda_i c[i]=0$.
Note that this implies in particular (since the function~$0$ is algebraic)
that the dimension of the linear span of
the functions $\{ \chi_c\, : \, c \in R\}$ equals the dimension of~$R$,
from which we deduce not only that the sum~$S \oplus Y$ is direct but also
that it equals the full solution space of the Mellin system.

In order to prove the above statement, we observe that if $\chi_\lambda$  is
an algebraic function of $(x_1, \dots, x_n)$ then, the homogenized function
$$\overline{\chi}_\lambda:=  \, \,\sum_{i=1}^s \lambda_i \sum_{k=1}^{m^{n-1}} \,
c[i]^{(k)}  \alpha^{(k)}$$
is an algebraic function of $a=(a_m, a_{m_1}, \dots,a_{m_n}, a_{n+1})$. 
But now, as in~\cite[Th.
3.5]{CDD},  we fix a point $a$ with $|a_{n+1}|$ sufficiently small
relative to the absolute values of the other coordinates and consider 
the loop~$\gamma$  defined by
$$ \gamma(\theta): = (a_m, a_{m_1}, \dots, a_{m_n}, {\rm exp}(2 \pi i m_n \theta) a_{n+1}).$$
Analytic continuation around~$\gamma$ will return all the local roots in a suitable
neighborhood of $a$ to their original values,
while  there is a subset~$S^{(k)}$ containing $m_n < m$ of the
roots, such that the logarithms of the roots in~$S^{(k)}$ will be
increased by~$2\pi i$. Then $\gamma^*(\alpha^{(k)})  -
\alpha^{(k)}$ equals~$2\pi i$ times the sum of the~$m_n$ roots
in~$S^{(k)}$, which cannot be zero since $R=R'$. We denote this sum of
roots by~$z_{S^{(k)}}$. 
Therefore
$$
\gamma^*(\overline{\chi}_\lambda) - \overline{\chi}_\lambda =
\gamma^*\left( \sum_{i=1}^s \lambda_i \sum_{k=1}^{m^{n-1}} c[i]^{(k)} \alpha^{(k)}\right)
- \sum_{i=1}^s \lambda_i  \sum_{k=1}^{m^{n-1}} c[i]^{(k)} \alpha^{(k)},
$$
which, after dividing by~$2 \pi i$, equals
$$
\sum_{i=1}^s \lambda_i \sum_{k=1}^{m^{n-1}} c[i]^{(k)} z_{S^{(k)}} =
\sum_{k=1}^{m^{n-1}}
\left(\sum_{i=1}^s \lambda_i c[i]^{(k)} \right) z_{S^{(k)}}.
$$
Since~$z_S^{(k)}$ is not the complete sum of all the roots corresponding
to one equation, the only way in which a finite iterate of the
action of~$\gamma$ returns~$\overline{\chi}_\lambda$ to its original value is
when $ \sum_{i=1}^s \lambda_i c[i]^{(k)} = 0$ for all~$k$, i.e.,
$\sum_{i=1}^s \lambda_i c[i]=0$, as asserted. \hfill $\square$

\vskip0.2cm

We also deduce:

\begin{corollary}
The monodromy representation of the Mellin system~(\ref{mellin})
is always reducible.
\label{monodromyCor}
\end{corollary}
\begin{proof}
Suppose that $n>1.$ By Theorem~\ref{structurethm}, the dimension of
the space of analytic solutions to~(\ref{mellin}) equals~$m^n.$
Since the solutions to the algebraic equation~(\ref{algebraic}) span
a linear space of dimension at most~$m,$ there is a nontrivial
invariant subspace in the solution space of~(\ref{mellin}).
If $n=1$ and $m_1 < m-1,$ then the solutions of the algebraic equation
form a nontrivial invariant subspace of dimension~$m-1.$
Finally, for $n=1$ and $m_1=m-1,$ the one-dimensional linear space
spanned by the rational solution~$x_1$ is invariant under the action
of the monodromy.
\end{proof}


\section{The case of dimension one}
\label{univariatesec}

In this section we further investigate the special case of the
Mellin system of equations associated with the trinomial algebraic
equation
\begin{equation}
y^m + x y^{m_1} - 1 = 0.
\label{trinomial}
\end{equation}
In this case the Mellin system~(\ref{mellin}) consists of a single
ordinary differential equation of order~$m,$ namely the equation
\begin{equation}
\left(
\prod_{k=0}^{m_1 -1}
(m_1 \theta + mk + 1)
\prod_{k=0}^{m - m_1 -1}
( (m - m_1) \theta + mk - 1) -
(-1)^{m_1} m^m \frac{d^m}{d x^{m}}
\right)
y(x) = 0,
\label{ordinaryMellin}
\end{equation}
where $\theta = x\frac{d}{dx}.$
Let $m = d a,$ $m_1 = d b,$ where $d={\rm GCD}(m,m_1).$
The discriminant of the algebraic equation~(\ref{trinomial}) (computed
with respect to~$y$) is given by $x^a b^b (a - b)^{a-b} - (-1)^b a^a.$
For $d=1,$ this expression coincides with the leading coefficient in
the differential equation~(\ref{ordinaryMellin}).
Throughout this section, we will call
the differential operator in the left-hand side of~(\ref{ordinaryMellin})
{\it the Mellin operator} and denote it by~$M(m,m_1).$
Special instances of this operator are considered in
Examples~\ref{ex21},\ref{ex32},\ref{ex31}, and~\ref{ex2kk}.

\begin{example}
\rm
Consider the cubic equation
\begin{equation}
y^3 + xy^2 - 1 = 0.
\label{eq32}
\end{equation}
The corresponding Mellin equation can be written in the form
\begin{equation}
(4x^3 - 27) y''' + 18 x^2 y'' + 4 xy' - 4y = 0.
\label{mellin32}
\end{equation}
A basis in the solution space of the differential
equation~(\ref{mellin32}) is given by the functions
$$
y_1(x) = x,  \qquad
y_2(x) = \sqrt[3]{108 - 8x^3 + 12\sqrt{81 - 12x^3}}, \qquad
y_3(x) = x^2/y_2(x),
$$
or by the roots of the algebraic equation~(\ref{eq32}).

Notice that after multiplication with~$x^2$ the differential
operator defining the equation~(\ref{mellin32}) can be
factorized as follows:
$$
x^2
\left(
(4x^3 - 27) \frac{d^3}{dx^3} + 18 x^2 \frac{d^2}{dx^2} + 4 x\frac{d}{dx} - 4
\right) = \phantom{------------}
$$
$$
\phantom{-----------}
\left(
(4x^4 - 27x)\frac{d^2}{dx^2} + (14x^3 + 27)\frac{d}{dx} + 4x^2
\right)
\left(
x\frac{d}{dx} - 1
\right).
$$
This factorization makes the rational solution $y_1(x)=x$ obvious.
\label{ex32}
\end{example}

More in general, assume $m_1 = m-1$.
Since the  roots $y_1(x), \ldots, y_m(x)$
of the trinomial satisfy
the equation~(\ref{ordinaryMellin}) and sum up to~$x,$
it follows that the operator~$\theta - 1$ (which is the annihilator
of~$x$) can be factored out of~$M(m,m-1).$ This factorization is
described in the following proposition.

\begin{proposition}
For $m_1 = m-1,$
the only nontrivial decomposition of the Mellin operator defining
the equation~(\ref{ordinaryMellin}) in the Weyl algebra of linear
differential operators with polynomial coefficients is given by
\begin{equation}
x^m M(m,1) =
\left(
x^m
\prod_{k=0}^{m - 2}
( (m - 1) \theta + mk + 1) + (-m)^m \, \theta \prod_{k=2}^{m-1}(\theta - k)
\right)(\theta - 1).
\label{m:m-1:Decomposition}
\end{equation}
\label{m:m-1:DecompositionProp}
\end{proposition}
The proof of this proposition is analogous to the proof of
Proposition~\ref{MellinDecompositionProp} below.

Assume now that $m > m_1 + 1$ and $ d= {\rm GCD}\,(m,m_1) = 1.$

\begin{example}
\rm
Consider the cubic equation
\begin{equation}
y^3 + xy - 1 = 0.
\label{eq31}
\end{equation}
The corresponding Mellin equation can be written in the form
\begin{equation}
(4x^3 + 27)y''' + 18x^2 y'' + 10xy' - 2y = 0.
\label{mellin31}
\end{equation}
Since the coefficient by~$y^2$ in~(\ref{eq31}) is zero, the
solutions to~(\ref{eq31}) generate a two-dimensional vector space
and thus there exists a solution~$\chi$ to the differential
equation~(\ref{mellin31}) which is not a solution to the algebraic
equation~(\ref{eq31}). The following functions give a basis of
solutions to the Mellin equation~(\ref{mellin31}):
$$
z_1(x) = \sqrt[3]{108 + 12\sqrt{12x^3 + 81}}, \qquad
z_2(x) = x/z_1(x),
$$
$$
\chi(x) = y_1 \log (y_{1}) + y_{2} \log ({y_2}) +  y_{3} \log({y_3}),
$$
where $y_1(x), y_2(x), y_3(x)$ are the three (linearly dependent)
solutions to~(\ref{eq31}).

Notice that the differential operator defining the
equation~(\ref{mellin31}) can be factorized as follows:
$$
(4x^3 + 27)\frac{d^3}{dx^3} + 18x^2 \frac{d^2}{dx^2} + 10x\frac{d}{dx} - 2 =
\frac{d}{dx}
\left(
(4x^3 + 27)\frac{d^2}{dx^2} + 6x^2 \frac{d}{dx} - 2x
\right).
$$
The functions~$z_1,z_2$ give a fundamental system of solutions to
the second-order factor in this decomposition. This operator is
described in~\cite{KatoNoumi}.
\label{ex31}
\end{example}

It turns out, that the Mellin operator corresponding to~$m_1 = 1$
allows a similar decomposition.

\begin{proposition}
For~$m_1 = 1,$
the only nontrivial decomposition of the Mellin operator defining
the equation~(\ref{ordinaryMellin}) in the Weyl algebra of linear
differential operators with polynomial coefficients is given by
\begin{equation}
M(m,1) =
\frac{d}{dx}
\left(
x \prod_{k=0}^{m - 2}
( (m - 1) \theta + mk - 1) + m^m \frac{d^{m-1}}{d x^{m-1}}
\right).
\label{MellinDecomposition}
\end{equation}
\label{MellinDecompositionProp}
\end{proposition}
\begin{proof}
The validity of~(\ref{MellinDecomposition}) follows from the Weyl
algebra identity $\theta + 1 = \frac{d}{dx} \circ x.$ Here~'$\circ$'
denotes the composition of differential operators in the Weyl algebra.
Let us denote the second factor in the right-hand side
of~(\ref{MellinDecomposition}) by~$M_{\rm alg}.$
We would like to show that the differential equation
$M_{\rm alg}y = 0$ has irreducible monodromy.
Multiplying this equation with~$x^{m-1},$
using the identity~(\ref{useful}) and making
the change of variables $t=x^m,$ we arrive at the equation
\begin{equation}
\left(
(m-1)^{m-1} \, t
\prod_{i=0}^{m-2}
\left(
\theta_t + \frac{mi-1}{m(m-1)}
\right) +
m^m
\prod_{j=0}^{m-2}
\left(
\theta_t - \frac{j}{m}
\right)
\right)y = 0,
\label{secondFactor}
\end{equation}
where $\theta_t = t\frac{d}{dt}. $
By Proposition~3.3 in~\cite{BeukersHeckman} the monodromy of
the differential equation~(\ref{secondFactor}) is reducible
if and only if
$(mi-1)/m(m-1) + j/m \in\Z,$ for some $i,j = 0,\ldots, m-2.$
This condition can only be satisfied when
$(mi-1)/m(m-1) + j/m = 1.$ Clearing the denominators and
considering the obtained relation modulo~$m$ we conclude that it
can never be satisfied. Thus the monodromy of~(\ref{secondFactor})
is irreducible and hence~$M_{\rm alg}$ does not admit any further
decomposition.
\end{proof}

In case $d > 1$ the Mellin operator always factors.

\begin{proposition}\label{singer}
There exist~$d$ linear operators $Q_0, \dots, Q_{d-1} $  of order~$\frac m d$
with coefficients in~$\C(x)$ such that $M(m,m_1) = Q_0 \dots Q_{d-1}$.
\end{proposition}

\begin{proof}
According to Corollary~\ref{gcd} and the discussion in Section~\ref{rootsec}, all
solutions to $M(m,m_1)$ are algebraic (so its associated Galois group coincides
with the monodromy group) and a basis of solutions is given by the
functions
$$
\{ f_{j,k}(x) :=\e^j y_{pr}(\e^{jm_1 + k} x), \, j=0, \dots, m/d-1, \, k=0, \dots,d-1\},
$$
where~$\e$ is a primitive $m$-root of~$1$ and~$y_{pr}$ is the principal solution
described in~(\ref{principalseries}).
For each $k =0, \dots, d-1$, the~$m/d$ functions $\{f_{j,k}(x), \, j=0, \dots,m_{d-1} \}$
are the (linearly  independent) roots of the algebraic equation
$$
y^m + \e^k x y^{m_1} -1 =0, \quad k=0, \dots, d-1.
$$
Then, they span a monodromy invariant subspace~$V_k$ of the
solutions~$V$ to~$M(m,m_1)$ of dimension~$m/d$, which gives a
right factor~$L_k$ of order~$m/d$ of the Mellin operator whose
solution space equals~$V_k$. Moreover, we have the direct sum
decomposition $V \, = \, V_0 \oplus \dots \oplus V_{d-1}$. Then,
it follows by the results in~\cite{Singer} that there exist
operators $Q_0, \dots, Q_{d-1}$ (with $Q_i$ equivalent to~$L_i$ in
the sense of~\cite{Singer}) such that~$M(m,m_1)$ decomposes as
$M(m,m_1) = Q_0 \dots Q_{d-1}$.
\end{proof}

\begin{example}
\rm
Let $n=1,$ $m=4,$ $m_1 = 2.$
The Mellin equation associated with the algebraic equation
$y^4 + x y^2 - 1 = 0$ is given by
\begin{equation}
(16x^4 - 256)y^{(4)} + 160 x^3 y''' + 360 x^2 y'' + 120 xy' - 15y =0.
\label{eq42}
\end{equation}
In accordance with Corollary~\ref{prsolcor}, a basis in the solution
space of~(\ref{eq42}) is given by solutions to the algebraic
equations $y^4 + x y^2 - 1 = 0$ and $y^4 + i x y^2 - 1 = 0.$
After a suitable normalization this basis is given by the functions
$\sqrt{-x \pm \sqrt{x^2 \pm 4}}.$ The solution space of~(\ref{eq42})
splits into two 2-dimensional invariant subspaces spanned by
$\sqrt{-x \pm \sqrt{x^2 + 4}}$ and $\sqrt{-x \pm \sqrt{x^2 - 4}}.$
This decomposition of the solution space corresponds to the
factorization of the differential operator defining~(\ref{eq42}):
$$
(16x^4 - 256)\frac{d^4}{dx^4} + 160 x^3 \frac{d^3}{dx^3} +
360 x^2 \frac{d^2}{dx^2} + 120 x \frac{d}{dx} - 15 =
\phantom{-------}
$$
$$
\phantom{-------}
\left(
(4x^2 - 16)\frac{d^2}{dx^2} + 20 x \frac{d}{dx} + 15
\right)
\left(
(4x^2 + 16)\frac{d^2}{dx^2} +  4 x \frac{d}{dx} - 1
\right).
$$
\label{ex42}
\end{example}
\begin{example}
\rm
Let $n=1,$ $m=6,$ $m_1 = 2$ and consider the algebraic equation
\begin{equation}
y^6 + xy^2 - 1 = 0.
\label{eq62}
\end{equation}
The Mellin equation associated with~(\ref{eq62}) is given by
$$
(1024 x^6 - 46656) y^{(6)} + 27648 x^5 y^{(5)} + 242816 x^4 y^{(4)} +
\phantom{---------}
$$
\begin{equation}
\phantom{----}
818944 x^3 y''' + 955780 x^2 y'' + 236180 x y' - 6545 y =0.
\label{mellin62}
\end{equation}
A basis in the space of solutions to~(\ref{mellin62})
is given by $\{ f(\varepsilon^i x), \, i=0,\ldots, 5 \},$
where $\varepsilon=e^{\pi i/3}$ and~$f(x)$ is the principal
solution to~(\ref{eq62}).
The differential operator defining the equation~(\ref{mellin62})
admits the following factorization:
$$
(1024 x^6 - 46656) \frac{d^6}{dx^6} + 27648 x^5 \frac{d^5}{dx^5}  +
242816 x^4 \frac{d^4}{dx^4} + \phantom{-----}
$$
$$
\phantom{-------}
818944 x^3 \frac{d^3}{dx^3} + 955780 x^2 \frac{d^2}{dx^2} +
236180 x \frac{d}{dx} - 6545 =
$$
$$
\left(
(32 x^3 - 216) \frac{d^3}{dx^3} + 432 x^2 \frac{d^2}{dx^2} +
1526 x \frac{d}{dx} + 1309
\right) \circ \phantom{---------}
$$
$$
\phantom{--------}
\left(
(32 x^3 + 216) \frac{d^3}{dx^3} + 144 x^2 \frac{d^2}{dx^2} +
86 x \frac{d}{dx} - 5
\right).
$$
\label{62}
\end{example}



\end{document}